\documentclass[12pt]{article}

\usepackage{amssymb}
\usepackage{amsmath}
\usepackage{amsthm}

\usepackage{tikz-cd}
\usepackage{comment}
\usepackage{carolmin}

\def\finb{{\rm Fin}(\beta)}

\def\N{\mathbb{N}}
\def\R{\mathbb{R}}
\def\Z{\mathbb{Z}}

\newtheorem{theorem}{Theorem}[section]
\newtheorem{corollary}[theorem]{Corollary}
\newtheorem{lemma}[theorem]{Lemma}
\newtheorem{proposition}[theorem]{Proposition}

\newtheorem{definition}[theorem]{Definition}
\newtheorem{remark}[theorem]{Remark}
\newtheorem{example}[theorem]{Example}

\begin{document}

\title{On finite $\beta$-expansions for the set of natural numbers}

\author{Túlio O. Carvalho\thanks{Corresponding author. E-mail: tuliocarvalho@uel.br} \and Catharina M. Moreira\\
Departamento de Matemática \\
Universidade Estadual de Londrina \\
Rod. Celso Garcia Cid, PR-445, Km 380 \\
CEP 86055-900, PR, Brazil}

\maketitle 
\begin{abstract}
We present a study of the problem of finiteness of the $\beta$-expansions for the set of natural numbers, condition $F_1$ in brief, for a cubic family of Pisot numbers for which the $\beta$-expansion of 1 is not a non-decreasing sequence. This is a class of simple $\beta$-numbers, within which a puzzling behaviour is displayed, in the sense that an infinite subset of such $\beta$ satisfy $F_1$, whereas a complementary infinite subset does not. This puzzle is organized by the residue class modulo 3 of $\lfloor \beta\rfloor$, which appears as a coefficient in the family of polynomials studied. We prove that,  when $F_1$ does not hold, the complement of ${\rm Fin}(\beta)$ is infinite. We eventually observe that, within this family, a concise sufficient condition for $F_1$ holds: the finitude of the $\beta$-expansion of $(\beta-\lfloor \beta\rfloor)^2$.
\end{abstract}






\section{Introduction}

Given a non-integer number $\beta>1$, 
Rényi \cite{renyi}, and subsequently Parry \cite{parry},
initiated the study of the dynamics of the piecewise expanding map of the interval: 
	\begin{align*}
		T_\beta:[0,1]&\rightarrow[0,1]\\
		x&\mapsto  \beta x -\lfloor\beta x \rfloor \ ,
	\end{align*}
where $\lfloor y\rfloor$ denotes the integer part of $y$. 
To each $x\in [0,1]$, write $e_1(x)=\lfloor \beta x\rfloor$ and $e_{k+1}(x)=\lfloor \beta T_\beta^k(x)\rfloor$, $k\in \N$. The integers $e_k(x)$ form a sequence which encode the {\em orbit} of $x$ under $T_\beta$. These are the digits, from the alphabet ${\cal A}=\{0,1,\ldots,\lfloor \beta\rfloor\}$, of the $\beta$-{\em expansion} of $x$, in the sense that
\begin{equation}
x = \sum_{k=1}^\infty \frac{e_k(x)}{\beta^k} \ . \label{greedy}
\end{equation}
Eq.\eqref{greedy} is also called the {\em greedy} expansion of $x$ in base $\beta$. 

When $x>1$, there is a least $L\in \N$ such that $x< \beta^L$, and the $\beta$-expansion of $x$ is defined as the $L$ times (left) shifted $\beta$-expansion of $\frac x{\beta^L}$.

Given a bounded sequence ${\bf c}=(c_1,c_2,\ldots)$ of integers, finite or not, we write
\[ {\bf c}(\beta)= \sum_{k=1}^\infty \frac{c_k}{\beta^k} \ . \]
Let us term a sequence {\em canonical} if ${\bf c}=(\lfloor \beta T^{k-1}(x)\rfloor)_{k\in \N}$ for some $x\in [0,1)$. 
A bounded sequence of integers ${\bf c}$ such that ${\bf c}(\beta)=x$ is called a $\beta$-{\em  representation} of $x$. 
Given $x>0$, we use the notation $[x]_{\beta e}$ for the possibly shifted $\beta$-expansion of $x$, and $[x]_\beta$ for any $\beta$-representation of $x$, including its $\beta$-expansion. 

Denote by ${\cal B}^\N$ the space of sequences in the alphabet ${\cal B}$.
Parry \cite{parry} gave a characterization of canonical sequences as a subshift of ${\cal A}^\N$. Assuming the usual order in ${\cal A}$, 
denote by $<_{\rm lex}$ the lexicographical order in the space of sequences ${\cal A}^\N$.
Introduce $\sigma(c_1,c_2,c_3,\ldots )=(c_2,c_3,\ldots)$ for the left-shift in ${\cal A}^\N$, and, given a finite sequence $t$, the notation $t^\omega$ to mean $ttt\ldots$. 
In the particular case where $\beta$ is a simple $\beta$-number, let
\begin{equation}
\label{betaexp1} 
[1]_{\beta e}=t_1t_2\ldots t_n \ .
\end{equation}
Then ${\bf c}$ is canonical if, and only if, for every $k\geq 0$, 
$$\sigma^k({\bf c})<_{\rm lex} (t_1t_2 \ldots (t_n-1))^\omega \ .$$

Denote the set of numbers having a finite $\beta$-expansion by $\finb$; $\Z[\beta^{-1}]$ the values of the ring of polynomials in $\beta^{-1}$ with integer coefficients; 
$\N[\beta^{-1}]$ the values of the positive cone in this ring and 
$\Z[\beta^{-1}]_+=\{x\geq 0\ :\ x\in \Z[\beta^{-1}]\}$.  
Following \cite{frougsolo}, we study the finiteness conditions 
\begin{align*}
\tag{$F_1$} & {\mathbb N} \subset \finb \ ,\\
\tag{$PF$} & \N[\beta^{-1}]  \subset \finb \ ,\\
\tag{$F$} & \Z[\beta^{-1}]_+  = \finb \ ,
\end{align*}
for Pisot numbers $\beta$ in some classes of polynomials. Frougny and Solomyak obtained an algorithm which proved condition $F$ for $\beta$ in a class of Pisot numbers studied by Brauer \cite{brauer50}. 

Akiyama \cite{akiyama} proved a criterion suitable for deciding whether condition $F$ holds for a particular Pisot number $\beta$: \medskip 

\noindent{\bf Theorem (Akiyama)}. {\it  
Let $\beta$ be a Pisot number. If every element of 
\[ {\cal C}= \left\{ x\in \Z[\beta]\ |\ 0<x<1,\; |x'|\leq 
\frac{\lfloor \beta\rfloor }{1-|\beta'|} \right\} \ , \]
where $\beta'$ runs through every Galois conjugate of $\beta$, has finite 
$\beta$-expansion, then $\beta$ has property $F$.} 
\medskip 

The studies concerned with the finiteness conditions have been guided by examples, first from Pisot numbers where $PF$ holds, but $F$ does not \cite{akiy_pf}. Akiyama stated then that ``{\em the full characterization of $\beta$ with $F$ among algebraic integers is a difficult problem}'', when the degree of the minimal polynomial of $\beta$ is greater than 2. 
It is now understood that $F \subset PF \subset F_1$, where the inclusions are all strict. While there are several contributions concerning cubic Pisot numbers with $PF$ without $F$, our paper is to provide examples where $F_1$ holds without positive finiteness $PF$. 

Sufficient conditions for $F_1$ were given in \cite{taka1}, using shift radix system (SRS). We will also come accross a sufficient condition for $F_1$, for the specific family of cubic polynomials we dealt with. This work may be viewed as a step towards generalizing extensive Akiyama's results on cubic Pisot units \cite{akiyama_cub}. 
Our approach to proving condition $F_1$ in section \ref{newver} will use a conjugation of the restriction of the map $T_\beta$ with an affine transformation in $\Z^3$, for a cubic polynomial. To be sure, the final argument coincides with the SRS type contraction argument.

Let $\beta$ be an algebraic integer root of the polynomial
\[ p(x)=x^m-a_1x^{m-1} -\cdots -a_{m-1}x-a_m \] 
and such that it is a simple $\beta$-number with the $\beta$-expansion of 1 given by eq.\eqref{betaexp1}.

One can consider various classes of algebraic numbers, defined by different conditions on the polynomial $p$ and/or $[1]_{\beta e}$: 
\begin{itemize}
\item[a)] If $p$ is irreducible and $|\beta-a_1|> 1+|a_2|+\cdots+|a_m|-|\beta|$ are Perron-Pisot numbers;
\item[b)] $a_1\geq a_2\geq \ldots \geq a_m>0$ are Brauer-Pisot numbers;
\item[c)] $a_1>a_2+\cdots+a_m$, $a_j\geq 0$, $\forall j$. 
\item[d)] $\beta$ is a Pisot number and the $\beta$-expansion of 1 is such that $(t_i)_{1\leq i \leq n}$ is non-monotonic and some $t_k=0$, $1<k<n$.
\end{itemize}
Perron \cite[p.~291, eq.5]{perron1907} was the first to consider condition (a) on Pisot numbers. As already mentioned, Brauer's condition (b) is sufficient for proving the finiteness condition $F$ \cite{frougsolo}. 
In class (b), $n=m$ and $t_i=a_i$ for each $1\leq i\leq n$. 
Condition (c) was proposed by Hollander \cite{hollander}, where $p$ is not necessarily irreducible.
Containing the counterexamples in \cite{frougsolo}, we name the class (d) {\em CE-Pisot numbers}. These are the Pisot numbers we will concern ourselves on studying $F_1$. Within the class CE, the smallest Pisot number with $\beta$-expansion $10001$ has property $F$. Though not defined, the class CE has already been studied in \cite{bassino,taka1}. In the latter paper, condition $F_1$, without $F$, is proved for the cubic polynomials: $x^3-2tx+2tx-t$, where $t$ is a natural number greater than 1.

For each natural number $N$, say its $\beta$-expansion $[N]_{\beta e}$ is the concatenation of a finite sequence $(d_{k-i})_{0\leq i\leq k}$ and an infinite sequence $(c_j)_{j\in \N}$, that is,
\begin{equation}
\label{suf_orb}
N= \sum_{i=0}^{k} d_i \beta^i+ \sum_{j=1}^\infty \frac{c_j}{\beta^j} \ , 
\end{equation}
If the sequence ${\bf c}=(c_j)_{j\in\N}$ has infinitely many non-zero terms, $N\not\in \finb$. We call ${\bf c}(\beta)$ the $\beta${\it -fractional part} of $N$. We will have occasion to show in passing that some subsets of natural numbers share the same $\beta$-fractional part.

The present study was induced by the remarkable counterexamples of \cite{frougsolo}: in the first $[1]_{\beta e}= 2102$ and $6$ has an infinite $\beta$-expansion, and in the second example $[1]_{\beta e}=2011$ and $3\not\in \finb$. The polynomials we use in our investigation do not belong in classes (a), (b) or (c), and we include some of degree four and five. 
At first, our interest was driven to find and show an infinite set of polynomials whose largest root would be a Pisot number $\beta$, simple in the sense of Parry, but for which condition $F_1$ would not hold. 

Within the family of counter-examples, there appears a subset of polynomials for which $F_1$ holds, a fact we prove by hand.
On the other hand, it seems remarkable that $F_1$ is true within the same family of polynomials where it is not. 
An important characteristic shared by the Pisot numbers $\beta$ whose $\beta$-expansions we study is that they are not units, which essentially implies $\beta^{-1}\not\in\Z[\beta]$, see Lemma \ref{beta_1}. 

The cubic family we studied is inspired by counterexamples from \cite{frougsolo}:
\begin{equation}
p_{n}(x) =x^3-(n+1)x^2+nx-n \ , \ n\geq 2\label{pol1} 
\end{equation}
As mentioned earlier, $p_2(x)= x^3-3x^2+2x-2$ has the leading root $\beta$, for which $[1]_\beta=2102$, but $6\not\in \finb$. 

\begin{theorem}
\label{zero}
Within the cubic family given in eq.\eqref{pol1}, there are infinitely many CE-Pisot numbers 
$\beta$ for which 
$[1]_{\beta e}$ is finite, but $F_1$ does not hold. 
\end{theorem}

The approach exposed in section \ref{newver} was developed mainly using the Pisot condition, which is similar (but different) to the approach using shift radix systems \cite{ak_all}. 

\begin{theorem}
\label{um}
Within the cubic family given in eq.\eqref{pol1}, there exist infinitely many CE-Pisot numbers $\beta$ for which $F_1$ holds, but $PF$ does not.
\end{theorem}

We stress the importance of the validity of $F_1$ coexisting, within the same family of polynomials, with its falsehood. However, the infinitude of cases of CE-Pisot numbers where $F_1$ is true is known since the 
recent work of \cite{taka1}. Our proof has a different reasoning. 

Usual sufficient conditions for proving that $F$ holds for a particular $\beta$, see \cite{akiyama,yosh1}, are decided upon investigation of finitely many $\beta$-expansions. One could be misled to think that once a natural number $k\not\in \finb$, then every natural number greater than $k$ also has an infinite $\beta$-expansion. However, when $F_1$ does not hold, we have the following remark.

\begin{remark}
\label{dois}
There exists infinitely many Pisot numbers $\beta$ such that $F_1$ does not hold, but the set ${\cal I}_\beta= \{k\in \N: k\not\in \finb\}$ has infinitely many gaps: 
$x\in {\cal I}_\beta$ is not a sufficient condition for $x+1\in {\cal I}_\beta$. 
\end{remark}

We will see that $\cal C$ in Akiyama's theorem include $n\beta+2\not\in \finb$, for any root $\beta$ in the cubic family $p_n$, preventing us from using this argument to settle Theorem \ref{um}.

Each polynomial $p_n$ is of Pisot type, that is, its root of largest module satisfies $\beta>1$, and the other roots lie inside the unit circle (see Lemma \ref{ispisot}). 
The $\beta$-expansions of 1 are finite, but do not conform with either Perron, Brauer or Hollander conditions.

We prove that condition $F_1$ is not true for $p_n$ when $n\not\equiv 1 \mod 3$ in section \ref{cubic}. 
The proof that it is true when $n\equiv 1\mod 3$ is concluded in section \ref{newver}.
As a corollary of Theorem \ref{um}, we make the following remark: 
\begin{remark}
\label{suff}
For a fixed $n$, let $\beta>1$ be the largest root of eq.\eqref{pol1} and $\alpha=\beta-\lfloor \beta\rfloor$. If $[1]_{\beta e}$ is finite, then $\alpha^2\in \finb$ is a sufficient condition for $F_1$. 
\end{remark}

The article is organized so that section \ref{cubic} contains proofs of Theorems \ref{zero} and Remark \ref{dois}, which are in a sense, negative results. The calculations are adaptations of the algorithm used to prove $\finb$ for Brauer-Pisot numbers explained in \cite{frougsolo}. Details are given paving the introduction of the restriction of the $\beta$-map to $\beta$-fractional parts of natural numbers: the content of section \ref{newver} necessary to prove Theorem \ref{um}. 
In section \ref{ultima}, a general construction based on the companion matrix of the cubic polynomial for $\beta$ is recalled to guide a discussion on other examples of Pisot numbers, outside the family given in eq.\eqref{pol1}.

\section{$\beta$-representations and unfolding words}
\label{cubic}

Calculations of the $\beta$-expansion of relatively small integers gives sufficient information for the disproof of condition $F_1$, and the content of Remark \ref{dois}. The $\beta$-expansion of $n^2+2$ reveals that for $n\equiv 0 \mod 3$ and $n\equiv 2 \mod 3$, condition $F_1$ does not hold. 

We follow \cite{frougsolo}. 
Given two finite sequences, $u=u_1\ldots u_m$, $v=v_1\ldots v_k$, padding zeros to the right of $u$ if necessary, we use the {\em positional summation} with the notation $\oplus_{j}$:
\[ u\oplus_j v = u_1\ldots u_{j-1}(u_j+v_1)(u_{j+1}+v_2) \ldots (u_{j+k-1}+v_k) u_{j+k}\ldots u_m \ . \]
Here $1\leq j\leq m$. We pad a zero to the left, that is $u_0=0$, when $j=0$.
In general, $u\oplus_j v$ is not the same as $v\oplus_j u$, since the position $j$ in $\oplus_j$ refers to the word on the left.

The positional summation satisfies the property:
\begin{align*}
	(u\oplus_j v) \oplus_k w= (u \oplus_k w)\oplus_j v \ ,
\end{align*}
By assuming the sequence of operations are done from left to right, we can drop the brackets. In seeking a $\beta$-expansion for a given finite $\beta$-representation, the positional summation $u\oplus_j v$ is used for sequences $v$ such that $v(\beta)=0$.

Given a natural number $N$, let $x$ denote its $\beta$-fractional part. 
Calculations aiming at its $\beta$-expansion depart from a known finite $\beta$-representation $[x]_\beta$.

Let us explain the algorithmic procedure, adapted from \cite{frougsolo}, supporting the calculations we  present. 
The first and main difference is in the domain of $\beta$-representations: we consider finite sequences $s$ in $\Z$, whose value $s(\beta)$ is positive. Since $s$ is a finite sequence, for some $k\in \N$, is ${\cal A}_k={\cal A}+\cdots +{\cal A}$ $k$ times, then $s\in -{\cal A}_k\cup {\cal A}_k$. 
If $u=u_1\ldots u_m$ is a finite sequence on some extended alphabet containing ${\cal A}$, its {\em length} is denoted by $|u|$, and equals $m$. The terminology in words follows reference \cite{lothaire}. 

If $|s|=k$, we write $s[1,k]=s$ and $s[i,j]=s_is_{i+1}\ldots s_j$, $1\leq i\leq j\leq k$. For instance, $\sigma(s)=s[2,k]$. Let $z=(-1)n10n$, which satisfies $z(\beta)=0$, for the roots of eq.\eqref{pol1}.

\noindent {\em Step 1.} 
Given $s$, first find the least index $j$ where $s[j,j+3]\geq_{\rm lex} z[2,5]$ (termwise). Search the first occurrence, say in position $i_1$, of a negative entry in $s[j+4,k]$. Find the first occurrence, say in position $i_2$, of a positive number in $s[j+4,k]$. Since $i_1\neq i_2$, we put
\begin{enumerate}
\item[a)] If $i_1<i_2$, we perform the operation $s\leftarrow s\oplus_{j+3} z$. 
\item[b)] Otherwise, let $s\leftarrow s\oplus_{j-1} (-z)$.
\end{enumerate}

Keep doing step 1, unless no word in $s$ is $\geq_{\rm lex} z[2,5]$ termwise.

\noindent {\em Step 2.}
Let $i_1$ be the least index such that $s[i_1]<0$, if it exists. Let $i_2$ be the least index such that $s[i_2,i_2+3]>_{\rm lex} z[2,5]$, if it exists. 

If there are no $i_1$ and $i_2$, the sequence $s$ is canonical and the recursion stops. 
\begin{enumerate}
	\item[a)] If $i_1<i_2$ or $i_2$ does not exist, we perform the operation $s\leftarrow s\oplus_{i_1-1} z$. 
	\item[b)] If $i_2<i_1$ or $i_1$ does not exist, there exists $0\leq p<3$ such that 
	\[ s[i_2,i_2+p]=z[2,2+p]\ \text{ and }\ s[i_2+p+1]> z[2+p+1] \ .\]
Set $s\leftarrow s\oplus_{i_2-1} (-z)\oplus_{i_2+p+1} z$.
\end{enumerate}

Given two words in an alphabet ${\cal A}$, $u=u_1\ldots u_m$ and $v=v_1\ldots v_k$, $v$ is a {\em suffix} of $u$ if 
\[ u_{m-k+1}\ldots u_m = v_1\ldots v_k \ . \]
$v$ is a {\em proper} suffix of $u$ if $m>k>0$. 
$v$ is a {\em prefix} of $u$ if $u_1\ldots u_k=v$, and $v$ is a {\em proper prefix} of $u$ if $|u|>|v|$.

Denote the length of the canonical prefix of $s$ by $c_s$. Let us examine what happens to $c_s\geq 0$ when submitted to this algorithm. After a set of runs of step 1a, $c_s$ will either increase or the sequence will enter step 1b. In step 1b, $c_s$ may decrease at first. The sum of digits decrease at each run of step 1b, so it cannot be performed indefinitely, since $s(\beta)>0$. 
After a set of runs of step 1b, $c_s$ increases. Analogously, on exiting a set of runs of step 2a, $c_s$ increases or the sequence enters step 2b. Finally, step 2b clearly increases $c_s$. The fundamental difference is that when $[1]_{\beta e}$ 
does not satisfy Brauer's condition, this algorithm can cycle, as the following example shows.

\begin{example}
\label{typ_cyc}
Take $s=n2$, and $[1]_{\beta e}=n10n$. In the first round of the procedure, step 2b, produces $s=1.00n(1-n)0n$.
The second round, $s$ enters step 2a, being modified to $s=1.00(n-1)11nn$. The third round (step 2b) produces
$s=1.00(n-1)120(n-2)n(1-n)0n$, and one sees the repetition of the sufix $n(1-n)0n$. 
\end{example}

\begin{definition}
Let $u$ and $v$ be two finite words. We say that $v$ is an {\em unfolding} of $u$, or that $v$ {\em unfolds} $u$, if $|v|>|u|$, and $u(\beta)=v(\beta)$. 
\end{definition}

For example, if $[1]_{\beta e}=2102$, the following calculation shows that $10102$ unfolds $11(-2)2$:
\begin{align*}
11(-2)2 & = 11(-2)2\oplus_{2}(-1)3(-2)2 \\
& = 10102 \ .
\end{align*}
Here, one has to note that since $\beta$ is a root of $p_2(x)=x^3-3x^2+2x-2$, then 
\[ 1 = \frac 3{\beta} -\frac{2}{\beta^2}+\frac{2}{\beta^3} \ , \]
so that $c=(-1)3(-2)2$ is such that $c(\beta)=0$.  

\begin{lemma}
\label{period_exp}
Suppose $v$ unfolds $u$, $u$ is a proper suffix of $v$: $v=wu$ for some non-empty word $w$, and $w$ is a canonical sequence satisfying $w<_{\rm lex} a_1\ldots a_{n-1}(a_n-1)$, where $a_1\ldots a_n = [1]_{\beta e}$. Then the $\beta$-expansion of $v(\beta)$ is infinite, periodic, and equals $w^\omega$. 
\end{lemma}

We may actually use that $\tilde{z}=(-1)(n+1)(-n)n$ also satisfies $\tilde{z}(\beta)=0$, when $\beta$ is a root of eq.\eqref{pol1}, 
and sum multiples of this sequence conveniently to improve the speed of the output of the algorithm, which is either a periodic or a finite $\beta$-expansion.

In order for $v$ to unfold $u$, we require $|v|>|u|$, which implies that unfolding is a transitive but not a reflexive relation. 
Let us prove that the maximal roots of eq.\eqref{pol1} are Pisot numbers. 

\begin{lemma}
	\label{ispisot}
	For each $n\in \N$, the maximal root of $p_n(x)$ given in eq.\eqref{pol1} is a Pisot number.
\end{lemma}

\begin{proof}
	$p_n(x)$ changes sign in the interval $(n,n+1)$. If $\beta$ denotes the largest positive root of $p_n$, as $p_n'(x)=3x^2-2(n+1)x+n>0$ for $x> n$, so that $n<\beta<n+1$. Let $x_1=\frac{n+1-\sqrt{n^2-n+1}}{3}$ denote the lesser root of $p_n'$, which is a local maximum for $p_n$. Then $0<x_1<1$ and $p_n(x)=-nx^2+(x^2+n)(x-1)$ enables us to see that $p_n(x_1)< -nx_1^2<0$. 
	
	This implies that $\beta$ is the unique real root of $p_n$. If $z$ is a complex root of $p_n$, then 
	\[ z\overline{z}\beta = n \ \Rightarrow\ |z|^2 = \frac{n}{\beta}<1 \ ,\]
	so that $\beta$ is a Pisot number.
\end{proof}

We analyse property $F_1$ validity dividing $\N$ in three intervals: 
\begin{enumerate}
\item[a)] $[2,n^2-n]$,
\item[b)] $[n^2+2,\infty)$ and 
\item[c)] $[n^2-n+1,n^2+1]$.
\end{enumerate} 
Every natural number in the set $[2,n^2-n]$ have finite $\beta$-expansion (Proposition \ref{pro0}). The analysis of the infinite set is inductive, see Proposition \ref{prep_induction}. The parametrization of the $\beta$-fractional parts which emerges extends back to $[n^2-n+1,n^2+1]$. 

Given a natural number $N$, we depart from one of its $\beta$-representations, and make positional sums of $\beta$-representations of zero, thereby increasing the length of its canonical prefix. The arrival on the unique $\beta$-expansion, when periodic, is guided by Lemma \ref{period_exp}. On the other hand, to show that $F_1$ holds, we will have knowledge about the possible $\beta$-fractional parts for every natural number $N$.

Recall the notation $\alpha=\beta-n$ for the $\beta$-fractional part of $\beta$. It turns out that the orbit of numbers of the form $1-k\alpha$, $k\in \N$ under the $\beta$-transformation is relevant for the $\beta$-expansion of natural numbers $2\leq N\leq n^2-n$. Note that $[\alpha]_\beta=0.10n$.

Given $k\in \N$, we have
\begin{align}
kn & = k(\beta-\alpha) = (k-1)\beta-(k-1)\alpha+n \nonumber \\
\label{lema1}
& = (k-1)\beta+(n-1)+(1-(k-1)\alpha) \ .
\end{align}
The significance of eq.\eqref{lema1} is in that it paves a way to deduce the $\beta$-expansions of natural numbers up to $n(n-1)$. 

Given $x>0$, in particular $x\in \N$, there is a finite word $u$ and a word $w$ such that $[x]_\beta = u.w$. We can assume that the prefix $u$ has been made canonical, after some positional summations if needed. The dot is to indicate that
\[ w(\beta)= \sum_{k=1}^\infty \frac{w_k}{\beta^k} \ . \]
If $0<w(\beta)<1$, $w$ is a $\beta$-representation of the $\beta$-fractional part of $x$. It is clear that determining whether $w(\beta)\in \finb$ is equivalent to establishing that $x\in \finb$. 

\begin{proposition}
\label{pro0}
For every $1\leq k\leq n(n-1)$, $k\in \finb$. 
\end{proposition}

The proof of Proposition \ref{pro0} is essentially contained in Lemma \ref{lema2} below. 

\begin{lemma}
If $1\leq k\leq n-1$, then $1-k\alpha \in (0,1)$. 
\label{obs1}
\end{lemma}

\begin{proof}
Since the statement follows from the definition of $\alpha$ for $k=1$, we just have to find a $\beta$-representation for $1-(n-1)\alpha$ with non-negative integers:
\begin{align*}
[1-(n-1)\alpha]_\beta & = 11(n-n^2)n \\
& = 11(n-n^2)n \oplus_{1} (-1)n10n \oplus_2 (-n)n^2n0n^2 \\
& = 01(n+1)(2n)0n^2 
\end{align*}
\end{proof}

Lemma \ref{obs1}, combined with \eqref{lema1}, shows that if $w_k$ is a $\beta$-representation of $1-(k-1)\alpha \in (0,1)$, the $\beta$-fractional part of $kn$, then 
\begin{equation}
\label{preobs}
[kn]_\beta = (k-1)(n-1).w_k \ . 
\end{equation}
Now Proposition \ref{pro0} is equivalent to proving that $w_k\in \finb$.

Now we examine the unfolding of a specific $\beta$-representation of $1-k\alpha$, for $1\leq k\leq n-1$:
\[ [1]_\beta=n10n,\; [k\alpha]_\beta= k0(kn)\; \Rightarrow\; [1-k\alpha]_\beta=
(n-k)1(-kn)n \ . \]
To establish Proposition \ref{pro0}, we must understand how these representations of $1-k\alpha$ evolve algorithmically. Lemma \ref{lema2} serves as the core unfolding mechanism for this purpose. It demonstrates that these words preserve their suffix structure when parametrized by integer numbers, while preserving canonical prefixes.

\begin{lemma}
$\forall j\in \N$, if $1\leq k\leq n-2$, $u_n(k,j)=(n-k)j(-kn)(jn)$ unfolds in $v_n(k,j)u_n(k,j+2)$, where 
\[ v_n(k,j)=(n-(k+1))(n-(k-j))(j+1)(k+1)(k-(j+1)) (n-(j+2)) , \]
and $v_n(k,j)$ is canonical for $1\leq k\leq n-2$, $1\leq j\leq k-1$.  
\label{lema2}
\end{lemma}

\begin{proof}
Note that $u_n(k,1)$ is the above referred $\beta$-representation of $1-k\alpha$ for $\beta$ Pisot number for the family \eqref{pol1}. 
We make a sequence of positional summations of zero: 
\medskip 

\noindent
\begin{tabular}{ccccccc}
$(n-k)$ & $j$ & $(-kn)$ & $jn$  \\
$-1$ & $n$ & $1$ & $0$ & $n$ \\
& $-k$ & $kn$ & $k$ & 0 & $kn$ \\
&& $j$ & $-jn$ & $-j$ & $0$ & $-jn$ \\
&&& $1$ & $-n-1$ & $n$ & $-n$ \\ \hline  
$(n-k-1)$ & $(n-k+j)$ & $(j+1)$ & $(k+1)$ & $(-j-1)$ & $(k+1)n$ & $(-j-1)n$
\end{tabular}
\medskip 

Continuing from the fifth position of the last word
\medskip 

\noindent 
\begin{tabular}{cccccc}
$(-j-1)$ & $(k+1)n$ & $(-j-1)n$ \\
$k$ & $-kn$ & $-k$ & $0$ & $-kn$ \\
& $-(j+2)$ & $(j+2)n$ & $(j+2)$ & $0$ & $(j+2)n$ \\ \hline 
$(k-j-1)$ & $(n-j-2)$ & $(n-k)$ & $(j+2)$ & $(-kn)$ & $(j+2)n$ 	
\end{tabular}
\medskip 

The conditions $1\leq k\leq n-2$ and $1\leq j\leq k-1$ provide for the six digits in $v_n(k,j)$ be non-negative. Moreover $v_n(k,j)$ and any of its five translations are $<_{\rm lex} n10n$.
\end{proof}

If $k$ is odd, the unfolding process of $u_n(k,1)$ eventually ends at the word with suffix $u_n(k,k)$, with the canonical prefix
\[ v_n(k,1)v_n(k,3)\ldots v_n(k,k-2) \ .\]
Now 
\begin{align*}
u_n(k,k)&= (n-k)k(-kn)(kn) \oplus_2 (-k) [k(n+1)](-kn)(kn) \\
&= (n-k)0k0(kn) 
\end{align*}
Recall that $u_n(k,1)$ originated from a $\beta$-representation of $1-k\alpha$, which by eq.\eqref{preobs} comes from a $\beta$-representation of the $\beta$-fractional part of $(k+1)n$. 

We have thus proved that if $k$ is odd and $kn\in \finb$, then $(k+1)n\in \finb$. 
If $k$ is even, we have to develop the suffix $u_n(k,k-1)$, which yields
\medskip 

\noindent 
\begin{tabular}{cccccccc}
	$(n-k)$ & $(k-1)$ & $-kn$ & $(k-1)n$ \\
	$-1$ & $n+1$ & $-n$ & $n$  \\
	& $-k$ & $kn+k$ & $-kn$ & $kn$  \\ 
	& $-1$ & $n$ & $1$ & $0$ & $n$ \\ 
	& & & $k$ & $-kn$ & $-k$ & $0$ & $-kn$ \\ \hline 
	$(n-k-1)$ & $(n-1)$ & $k$ & $(k+1)$ & $0$ & $(n-k)$ & $0$ & $-kn$ 	\\
\end{tabular}
\medskip 

Another sequence of positional summations (on the four letter suffix) leads to the word equation
\begin{align*}
0(n-k)0(-kn) = 0(n-k-1)(n-k)1(k+1)(k-1)(n-2)(n-k)2(-nk)(2n)
\end{align*}
Using Lemma \ref{lema2}, with $j=2$,
\begin{align*}
& (n-k)2(-nk)(2n)  = (n-k-1)(n-k+2)3(k+1) (k-3)(n-4)(n-k)
\ldots \\
& \hspace{2in} \ldots 
(n-k)k(-nk)(nk) \\
& \hspace{0.3in} = (n-k-1)(n-k+2)3(k+1) (k-3)(n-4)(n-k)
\ldots (n-k)0k0(kn)
\end{align*}
This proves that if $k$ is even and $kn\in \finb$, then $(k+1)n\in \finb$. 

To conclude the proof of Proposition \ref{pro0}, note that we showed that the $\beta$-expansions of $n,2n,\ldots , (n-1)n$ are all finite. If $1\leq j<n$, we have 
\[ n=(n-1).n10n \; \Rightarrow\; j= (j-1).n10n \ . \]
Analogously, applying  \eqref{preobs}, for $1\leq j<n$
\[ [(k+1)n]_\beta = k(n-1).w \; \Rightarrow\; [kn+j]_\beta = k(j-1).w\]
for some finite word $w$, show that natural numbers $N$ satisfying 
$kn+1\leq N\leq (k+1)n$ have the same $\beta$-fractional part. \qed

The arguments for the proofs of Theorem \ref{zero} and Remark \ref{dois} are built from moving forward the scrutiny of $\beta$-expansions of natural numbers around $n^2$.
One of its $\beta$-representations follows from equation \eqref{lema1}:
\[ n^2 = (n-1)\beta+(n-1)+ (1-(n-1)\alpha) \ . \]
The $\beta$-fractional part of $n^2$ is $(1-(n-1)\alpha)$, which has $11(n-n^2)n$ as a $\beta$-representation. One can see that $w(\beta)>0$ and, for $n\geq 3$, its $\beta$-expansion begins with 0. Indeed,
\begin{align*}
& \frac 1{\beta}+\frac 1{\beta^2}+\frac{n-n^2}{\beta^3}+\frac{n}{\beta^4} < \frac 1{\beta} \Leftrightarrow \beta+n+\frac{n}{\beta}< n^2 \ ,\; \text{ and} \\
& \beta+n+\frac{n}{\beta}< 2n+2< n^2
\end{align*}
It follows that if $n^2\in \finb$ so are the numbers $n(n-1)+1\leq N \leq n^2+1$, all of which have the same $\beta$-fractional part, for a fixed $n\geq 3$.

Since $n^2+2<\beta^2<n^2+3$, the operation of carry will transform the word to the left of the dot in a way that is different from the previous. This is due to forbidden words \cite{ambroz} appearing around the dot in the $\beta$-representations of these numbers obtained by summing 1 digitwise, that is, summing 1 to the digit multiplying $\beta^0$. For instance,  
\[ [n^2]_\beta = (n-1)(n-1).11(n-n^2)n \; \Longrightarrow\; 
[n^2+2]_\beta = n1.01(-n^2)n \ .\]
But since $1-\frac{n^2}{\beta}+\frac{n}{\beta^2}<0$, the $\beta$-fractional part of $n^2+2$ is given by
the word with value $2-n\alpha$: $n2(-n^2)(2n)$.

We depart from the fact that the $\beta$-expansion of any $N\in \N$ has a finite sequence to the left of the dot that separates the non-negative powers of $\beta$ from the negative ones (see eq.\eqref{suf_orb}). We state the following Proposition concerning the $\beta$-fractional parts of large enough natural numbers, which will allows us to think inductively on studying the behaviour of the $\beta$-fractional part of natural numbers larger than $n^2+1$.

\begin{proposition}
\label{prep_induction}
Let $p.w$ denote a $\beta$-representation of a number $N\geq n^2+2$, where $p=p_k\ldots p_0$ is canonical, and $0<w(\beta)<1$ is the $\beta$-fractional part of $N$. Let $(u_j)_{j\in \N}$ be a canonical word such that $u(\beta)=w(\beta)$. Suppose $p_k\ldots (p_0+1)$ is a forbidden word. Then, the following are the rules for calculating one 
$\beta$-representation, denoted by $w'$, of the $\beta$-fractional part of $N+1$:
\begin{itemize}
\item[a)] If $p_0=n$, then $w'=w\oplus_0(n-1)1(-n)n$.
\item[b)] If $p_0=n-1$, $p_3p_2p_1<_{\rm lex} n10$ and $u_1u_2u_3 \geq_{\rm lex} 10n$, then $w'=w\oplus_0(-1)0(-n)$. 
\item[c)] If $p_1p_0=n0$ and $u_1u_2\geq_{\rm lex} 0n$, then 
$w'=w\oplus_0 0(-n)=w\oplus_0(-1)1(-n)(n)$. 
\item[c')] If $p_1p_0=n1$, then $w'=w\oplus_0 0(-n)=w\oplus_0 (-1)1(-n)n$. 
\item[d)] If $p_2p_1p_0=n10$, then $w'=w\oplus_0 010n$. 
\item[e)] If $p_3p_2p_1p_0=n10(n-1)$, then $w'=w$.
\end{itemize}
\end{proposition}

\begin{proof}
The condition $N\geq n^2+2$ implies that the length of the canonical word $p_k\ldots p_0$ is at least 2.

In case (a), necessarily $w(\beta)<\alpha$, since otherwise case (b) would have been applied from $N-1$ to $N$. The digitwise operation in $w$ corresponds to   $w'(\beta)=w(\beta)+1-\alpha$.

Similarly, in case (d), we have $u_1<n$. Then $w'=w\oplus_0 n10n-n000$. 
\end{proof}

The effect of summing 1 to go from $n^2+1$ to $n^2+2$ illustrates case (a). The $\beta$-fractional part of $n^2+2$ is $n2(-n^2)n$. From $n^2+2$ to $n^2+3$, one uses case (c), reaching $(n-1)3(-n-n^2)(3n)$.

\begin{corollary}
For each $N\geq n^2+2$, its $\beta$-fractional part has a $\beta$-representation of the form
\[ \omega_n(l,k,j)= (ln-k)(j+l)(-kn-n^2)[(j+l)n] \ , \]
with $l\in \N$, $k\in \N$ and $j\in \Z$. 
\label{omega}
\end{corollary}

\begin{proof}
The $\beta$-fractional part of $n^2+2$ has $\beta$-representation 
$n2(-n^2)(2n)=\omega_n(1,0,1)$. 
If $\omega_n(l,k,j)$ denotes the $\beta$-fractional part of $N$, Proposition \ref{prep_induction} says the $\beta$-fractional part of $N+1$ is
\begin{itemize}
\item[(a)] $\omega_n(l+1,k+1,j)$,
\item[(b)] $\omega_n(l,k,j-1)$,
\item[(c,c')] $\omega_n(l,k+1,j+1)$,
\item[(d)] $\omega_n(l,k,j+1)$,
\item[(e)] $\omega_n(l,k,j)$. 
\end{itemize} 
\end{proof}

With a slight abuse of notation, we write $\omega_n(l,k,j)\in \finb$ when this four-lettered word represents a number with a finite $\beta$-expansion. 

We mention that Proposition \ref{prep_induction} is an illustration of the approach which lead Takamizo to Corollary 2.6 in \cite{taka1}.
However, this example gives a new ingredient with the appearance of modulo 4 arithmetics. Let us note that the parametrization of the $\beta$-fractional part extends back to $[n^2]_\beta$, since 
\[ \omega_n(0,-1,1)=11(n-n^2)n \ . \]

We have the following unfolding for $\omega_n(0,-1,j) = 1j(n-n^2)(jn)$:
\medskip

		\begin{tabular}{cccccccc}
			$1$ & $j$ & $n-n^2$ & $jn$ \\
			$-1$ & $n+1$ & $-n$ & $n$ \\
			& $-n+1$ & $n^2-n$ & $n-1$ & $0$ & $n^2-n$ \\ 
			& $-1$ & $n$ & $1$ & $0$ & $n$\\
			 & & $j+1$ & $-(j+1)n$ & $-(j+1)$ & 0 & $-(j+1)n$ \\ \hline 
$0$ & $j+1$ & $j+1$ & $n$ & $-(j+1)$ & $n^2$ & $-(j+1)n$
		\end{tabular} 
		\medskip 
		
		Continuing from position 5 of the last word:
		\medskip 
		
		\begin{tabular}{cccccc}
			$-(j+1)$ & $n^2$ & $-(j+1)n$ \\
			$n$ & $-n^2$ & $-n$ & $0$ & $-n^2$ \\
			$-1$ & $n+1$ & $-n$ & $n$ \\ 
			& $-(j+3)$ & $(j+3)n$ & $j+3$ & $0$ & $(j+3)n$ \\ 
			&& $1$ & $-n-1$ & $n$ & $-n$ \\ \hline 
			$n-j-2$ & $n-j-2$ & $1$ & $ (j+2)$ & $(n-n^2)$ & $(j+2)n$  
		\end{tabular}
		\medskip 
		
This provides the unfolding 
\begin{align*}
& \omega_n(0,-1,j)  = 1j(n-n^2)(jn)  \nonumber \\
&\; = 0(j+1)(j+1)n(n-j-2)(n-j-2)1(j+2)(n-n^2)[(j+2)n] 
\end{align*} 
If $j+3\leq n$, or $n-j-2\geq 1$, we get a canonical prefix by performing $\oplus_3 1(-n)(-1)0(-n)\oplus_6 (-1) n10n$ on the last $\beta$-representation of $\omega_n(0,-1,j)$: 
\begin{equation}
\label{unfold_n2}
\omega_n(0,-1,j) = 0(j+1)(j+2)0(n-j-3)(n-j-3)\omega_n(0,-1,j+3) \ .
\end{equation}
When $0\leq j\leq n-2$, the prefix $0(j+1)(j+2)0(n-j-3)(n-j-3)$ is canonical. 

\begin{remark}
	\label{alpha2} 
	The number $\alpha^2$, which has $10(n-n^2)$ as one of its $\beta$-representations, is parametrized as $\omega_n(0,-1,0)=[\alpha^2]_\beta$. 
	
	Application of eq.\eqref{unfold_n2} with $j=0$ gives a canonical prefix for its unfolding. 
\end{remark}

	Theorem \ref{zero} is a corollary of the first case in Proposition \ref{pro1}. 

\begin{proposition}
	\label{pro1}
	Let $\beta>1$ be the Pisot number of the polynomial $p_n$, $n\geq 2$, given in equation \eqref{pol1}. Then
	\begin{itemize}
		\item $n\equiv 0 \mod 3$ implies $n^2\not\in \finb$.
		\item $n\not\equiv 0 \mod 3$ implies $n^2 \in \finb$. 
	\end{itemize}
\end{proposition}

\begin{proof}
	The $\beta$-fractional part of $n^2$ is $\omega_n(0,-1,1)=11(n-n^2)n$. From the unfolding eq.\eqref{unfold_n2}, it unfolds with suffix $\omega_n(0,-1,1+3p)$, $p=0,1,\ldots$, and its prefix is canonical if $1+3p\leq n-2$. 
	The last suffix in this unfolding depends on the residue class modulo 3 of $n$. 
	
If $n\equiv 0 \mod 3$, $\omega_n(0,-1,1)$ unfolds till it reaches the suffix 
	$\omega_n(0,-1,n-2)= 1(n-2)(n-n^2)(n^2-2n)$. We can deduce that
	\begin{align*} 
	& 1(n-2)(n-n^2)(n^2-2n) = 0(n-1)(n-1)n010(n-n^2) \\
	&\quad = 0(n-1)(n-1)n 0\, \omega_n(0,-1,0) \ . 
	\end{align*}
	By the same rule eq.\eqref{unfold_n2}, $\omega_n(0,-1,0)$ unfolds till it reaches the suffix $\omega_n(0,-1,n)= 1 n (n-n^2)n^2$, actually preceeded by a $0$. 
	Now we note that $0\omega_n(0,-1,n)=01n(n-n^2)n^2 = 10(n-n^2)$:
	\begin{align*}
		& 10(n-n^2) \oplus_1 (-1)(n+1)(-n)n = \\
		& 0(n+1)(-n^2)n\oplus_{2} (-n)(n^2+n)(-n^2)(n^2) = \\
		& 0 1 n(n-n^2)n^2 \ .
	\end{align*}
	Using Lemma \ref{period_exp}, we conclude  that for $n\equiv 0\mod 3$, $n^2\not\in \finb$. 
	
	Consider now $n\equiv 1 \mod 3$, $\omega_n(0,-1,1)$ unfolds till it reaches the suffix $\omega_n(0,-1,n)=1n(n-n^2)n^2$, which we can show to unfold as
	\[ 1n(n-n^2)n^2= 1 2 0(n-3)(n-3) 13(n-n^2) (3n) = 120(n-3)(n-3) \omega_n(0,-1,3) \ . \]
	$\omega_n(0,-1,3)$ unfolds with the suffix $\omega_n(0,-1,n-1)$, which belongs to $\finb$, since
	\[ 1(n-1)(n-n^2)(n^2-n)\oplus_2 (1-n)(n^2-1)(n-n^2)(n^2-n) =
	10(n-1)0(n^2-n) \ . \]
	
	Finally, for $n\equiv 2 \mod 3$, $\omega_n(0,-1,1)$ unfolds till it reaches the suffix $\omega_n(0,-1,n-1)$, which we just noticed to belong to $\finb$. 
\end{proof}

Let us note that all the numbers $N$ in the interval $n^2-n+1\leq N\leq n^2+1$, for $n\equiv 0 \mod 3$ fixed, have the same  periodic suffix given by the word
\[ \big(0120(n-3)(n-3)0450(n-6)(n-6)\ldots 0(n-2)(n-1)00\big)^\omega \ , \]
by the closing argument in the proof of Proposition \ref{pro0}. On the other hand, this argument proves that, when $n\not\equiv 0 \mod 3$, the $\beta$-expansions of the natural number in the interval $[n^2-n+1,n^2+1]$ are finite.

We move now to examining the $\beta$-expansion of $n^2+2$. As mentioned earlier, its $\beta$-fractional part is 
$n2(-n^2)(2n)=\omega_n(1,0,1)$. The unfolding of $\omega_n(1,0,1)$ depends of the residue class modulo 4 of $n$, because with calculations similar to the unfolding of $\omega_n(0,-1,j)$ we get
\[ \omega_n(1,0,j-1) = (n-1)j(j+2)1(n-j-2)(n-j-4)\omega_n(1,0,j+3) \ . \]
The prefix is canonical when $0\leq j\leq n-4$. 
\medskip 

	\begin{tabular}{cccccccc}
	$n$ & $j$ & $-n^2$ & $jn$ \\
	$-1$ & $n+1$ & $-n$ & $n$ \\
	& $-n-1$ & $n^2+n$ & $n+1$ & $0$ & $n^2+n$ \\ 
	& & $j+2$ & $-(j+2)n$ & $-(j+2)$ & 0 & $-(j+2)n$ \\ \hline 
	$n-1$ & $j$ & $j+2$ & $1$ & $-(j+2)$ & $n^2+n$ & $-(j+2)n$
\end{tabular} 
\medskip 

Continuing from position 5 of the last word:
\medskip 

\begin{tabular}{cccccc}
	$-(j+2)$ & $n^2+n$ & $-(j+2)n$ \\
	$n$ & $-n^2-n$ & $n^2$ & $-n^2$ \\
	& $n-(j+4)$ & $-n^2+(j+4)n$ & $-n+j+4$ & $0$ & $-n^2+(j+4)n$ \\ 
	& & $-n$ & $n^2+n$ & $-n^2$ & $n^2$  \\ \hline 
	$n-j-2$ & $n-j-4$ & $n$ & $ (j+4)$ & $-n^2$ & $(j+4)n$  
\end{tabular}
\medskip 

This unfolding leads to four 
different possible behaviours for the suffix of $\omega_n(1,0,1)$ according to which residue modulo 4 of $n$. 

\begin{proposition}
\label{pro2}
According to the residue modulo 4 of $n$, 
the suffixes of $\omega_n(1,0,1)$ are
\begin{align*}
&n\equiv 0 \mod 4 \ , & & \omega_n(1,0,n-3)= n(n-2)(-n^2)(n^2-2n) \ ,\\
&n\equiv 1 \mod 4\ , & & \omega_n(1,0,n-4)=n(n-3)(-n^2)(n^2-3n) \ , \\
& n \equiv 2 \mod 4 \ , & &  \omega_n(1,0,n-1)=nn(-n^2)(n^2) \ , \\
&n \equiv 3 \mod 4\ , & & \omega_n(1,0,n-2)=n(n-1)(-n^2)(n^2-n) \ .
\end{align*}
\end{proposition}

\begin{proof}[Proof of Remark \ref{dois}]
Recall that $n^2\not\in \finb$, when $n\equiv 0 \mod 3$. 
The positional summation \medskip 

\begin{tabular}{ccccccccc}
  $n$ & $n$ & $-n^2$ & $n^2$  \\ 
   & $-n$& $n^2$ & $n$ & $0$ & $n^2$ \\ 
  & & $n$ & $-n^2-n$ & $n^2$ & $-n^2$ \\ \hline 
  $n$ & $0$ & $n$ & $0$ & $n^2$ \\
  & $1$ & $-n$ & $-1$ & $0$ & $-n$ \\
 &&& $n$ & $-n^2$ & $-n$ & $0$ & $-n^2$ \\ \hline
  $n$ & $1$ & $0$ & $n-1$ & $0$ & $-2n$& $0$& $-n^2$  \\ \hline
 &&& $-1$ & $n$ & $1$ & $0$ & $n$ \\
 &&&& $-2$ & $2n$ & $2$ & $0$ & $2n$ \\ \hline
  $n$ & $1$ & $0$ & $n-2$ & $n-2$ & $1$ & $2$ & $n-n^2$ & $2n$ 
\end{tabular}
\medskip

\noindent shows that the 4 lettered 
suffix of $\omega_n(1,0,1)$ when $n\equiv 2 \mod 4$ unfolds with suffix 
\[ 12(n-n^2)(2n) = \omega_n(0,-1,2) \ .\]

Applying eq.\eqref{unfold_n2}, for $j=2\equiv (-1)\mod 3$, this word will unfold till its suffix reaches the word $1(n-1)(n-n^2)(n^2-n)$. 
The next unfolding concludes the proof that $n^2+2\in \finb$:
\[ 
1(n-1)(n-n^2)(n^2-n) \oplus_2 (1-n)(n^2-1)(n-n^2)(n^2-n) = 10(n-1)0(n^2-n)
\ ,\]
since $n^2-n\in \finb$.

This behaviour reappears on the numbers which are in the pre-orbit of $n^2$ with respect to $\tau_\beta$, a map we will define in section \ref{newver}.
\end{proof}

We remark that $\omega_n(1,0,n-4)=n(n-3)(-n^2)(n^2-3n)$, the suffix for $n\equiv 1 \mod 4$, unfolds with suffix 
$\omega_n(0,-1,2)$ as well. When $n\equiv 0\mod 3$, $\omega_n(0,-1,2)$ has the suffix $1(n-1)(n-n^2)(n^2-n)$, which proves that $n^2+2$ also belongs in $\finb$ in this case. Thus $n=6$ and $n=9$ are the first examples in the family $p_n$ for Remark \ref{dois}.

\begin{proposition}
\label{n2p2}
If $n\equiv 2 \mod 3$ and $n\not\equiv 3\mod 4$, then $n^2+2\not\in \finb$. 
\end{proposition}

\begin{proof}
From the proof of Proposition \ref{pro2}, we know that for $n\equiv 2 \mod 4$, 
the four lettered suffix $nn(-n^2)(n^2)=\omega_n(1,0,n-1)$ unfolds in $\omega_n(0,-1,2)$. For $n\equiv 2 \mod 3$, from eq.\eqref{unfold_n2}, it unfolds till it reaches the suffix $\omega_n(0,-1,n)=1n(n-n^2)(n^2)$. The unfolding we just saw in the proof of Remark \ref{dois} shows that the suffix of $\omega_n(0,-1,n)$ goes back to $\omega_n(0,-1,2)$. This settles the proof for $n\equiv 2\mod 4$, by Lemma \ref{period_exp}.  

When $n\equiv 1 \mod 4$, $\omega_n(1,0,n-4)$, which we remarked that unfolds with suffix $\omega_n(0,-1,2)$. Since $n\equiv 2\mod 3$, it unfolds with suffix $\omega_n(0,-1,n)$ by eq.\eqref{unfold_n2}, which unfolds back to $\omega_n(0,-1,2)$. This settles the proof for $n\equiv 1\mod 4$, by Lemma \ref{period_exp}.

For $n\equiv 0 \mod 4$, we unfold the suffix $\omega_n(1,0,n-3)$:
\medskip

\begin{tabular}{cccccccc}
$n$ & $n-2$ & $-n^2$ & $n^2-2n$ \\
$-1$ & $n+1$ & $-n$ & $n$ \\
& $-n-1$ & $n^2+n$ & $n+1$ & $0$ & $n^2+n$ \\ \hline
$n-1$ & $n-2$ & $0$ & $n^2+1$ & $0$ & $n^2+n$ \\
&& $n$ & $-n^2$ & $-n$ & $0$ & $-n^2$ \\
&&&& $n$ & $-n^2-n$ & $n^2$ & $-n^2$ \\ \hline
$n-1$ & $n-2$ & $n$ & $1$ & $0$ & $0$ & $0$ & $-n^2$ \\
&&& $-1$ & $n$ & $1$ & $0$ & $n$ \\ \hline
$n-1$ & $n-2$ & $n$ & $0$ & $n$ & $1$ & $0$ & $n-n^2$ 
\end{tabular}
\medskip 

We recognize the suffix $\omega_n(0,-1,0)=[\alpha^2]_\beta$, whose $\beta$-expansion begins with $0$, but does not belong to 
$\finb$ when $n\equiv 2\mod 3$. 

However, for $n\equiv 3 \mod 4$, the unfolding of $\omega_n(1,0,n-2)$ is short:
\medskip 

\begin{tabular}{ccccc}
$n$ & $n-1$ & $-n^2$ & $n^2-n$ \\
$-1$ & $n+1$ & $-n$ & $n$ \\
& $-n$ & $n^2+n$ & $-n^2$ & $n^2$\\ \hline 
$n-1$ & $n$ & $0$ & $0$ & $n^2$
\end{tabular}
\medskip 

The $\beta$-expansion of $n^2$ is finite and the two consecutive zeros above will generate a canonical finite $\beta$-expansion for $\omega_n(1,0,n-2)$. 
\end{proof}

Under the same hypotheses of Proposition \ref{n2p2} and under similar calculations we can prove that $n^2+3\in \finb$. The $\beta$-fractional part of $n^2+3$ is $\omega_n(1,1,3)$.

\begin{remark}
\label{comple_prop214}
The first natural number whose $\beta$-expansion is infinite for $n\equiv 2 \mod 3$ and $n\equiv 3 \mod 4$ is $n^2+2n+3$. Its $\beta$-fractional part is $\omega_n(1,3,2)=(n-3)3(-3n-n^2)(3n)$. It has this sequence of suffixes: $\omega_n(1,3,n-1)$, then $\omega_n(1,1,0)$, then $\omega_n(1,1,n-3)$, then $n10(n-n^2)$. As in the proof of Proposition \ref{n2p2}, we recognize $[\alpha^2]_\beta = 10(n-n^2)$ preceded by $n$, whose $\beta$-expansion is infinite for $n\equiv 2\mod 3$.
\end{remark}

\section{A conjugate map of the $\beta$-map in $\Z[\beta]$}
\label{newver} 

The combinatorial framework based on positional summations and word unfolding, as developed in Section \ref{cubic}, provides an effective algorithmic tool for tracking the evolution of $\beta$-fractional parts of natural numbers and establishing that condition $(F_1)$ does not hold for $n\not\equiv 1 \pmod 3$. However, to prove that $(F_1)$ holds for $n \equiv 1 \pmod 3$, as stated in Theorem \ref{um}, a global analysis of the underlying dynamical system is required. To achieve this, we now transition from algorithmic string manipulations to an algebraic and geometric setting. Motivated by the parametrization of the $\beta$-fractional parts obtained in Corollary \ref{omega}, we construct an affine transformation on $\mathbb{Z}^3$ that acts as a conjugate of a restriction of the $\beta$-map. This conjugation enables us to investigate the periodic orbits and invariant subsets of the system, ultimately yielding contraction arguments necessary to complete the proof.

Let us introduce the functional $f:\Z^3\to \R$ given by
\begin{equation}
f(l,k,j) = l-\alpha k+ \frac{\alpha}{\beta}j \ , 
\label{func_cond}
\end{equation}
which is identified with the row vector of its coefficients $[1\ -\alpha\ \frac{\alpha}{\beta}]$. 
Note that, since $n=\beta-\alpha$ and $\beta$ is a root of $p_n(x)=x^3-(n+1)x^2+nx-n$, we get 
\begin{align*}
& f(1,n,-1) = 1-n\alpha -\frac{\alpha}{\beta} = \alpha^2-\alpha = -\frac{n^2}{\beta^3}
\end{align*}
Therefore 
\[ [f(l+1,k+n,j-1)]_\beta = \omega_n(l,k,j) \ . \]

Let $G$ be the $3\times 3$ matrix
\[ G= \begin{bmatrix} n & -1 & 0 \\ -1 & 1 & -1 \\ 0 & 1 & 0 \end{bmatrix} \]
representing the operator $G(l,k,j)=(nl-k,-l+k-j,k)$ from $\Z^3$ into itself.
The characteristic polynomial of $G$ is $p_n(x)$. Furthermore 
$$fG=[1\ -\alpha \;\ \frac{\alpha}{\beta}] G = \beta f \ ,$$ 
so that $f$ is a left-eigenvector of $G$ with eigenvalue $\beta$. 
We get a conjugate map of the $\beta$-map:
\begin{align*}
\beta f(l,k,j) & = f\big(G(l,k,j)\big) \ \Rightarrow\ \\
T_\beta(f(l,k,j)) &= f\big(G(l,k,j)\big) -\lfloor  f\big(G(l,k,j)\big)
\rfloor 
\end{align*}
so that $T_\beta$ acts as the shift on $\beta$-expansions the numbers $f(l,k,j)$, where $(l,k,j)$ belongs to 
\begin{equation}
U  =\{ (l,k,j)\in \Z^3 \ |\ f(l,k,j) \in [0,1)\} \ .
\label{pontosdez3}
\end{equation}
Denoting 
\begin{align*} 
\pi : &\;  \Z^3 \;\ \longrightarrow\;\  \R \\
& (l,k,j)  \mapsto  f(l+1,k+n,j-1)=\omega_n(l,k,j)(\beta) \ ,
\end{align*}

We note that $1,\, \alpha$ and $\alpha/\beta$ are rationally independent, since otherwise $\beta$ would have a minimal polynomial of degree 2 over the integers. Thence $\pi$ is injective. Also, denote by $U^*=U\setminus\{(0,0,0)\}$. 
We get the following commutative diagram 
\[ \begin{tikzcd}
\pi(U) \arrow{r}{T_\beta} & \pi(U) \arrow{d}{\pi^{-1}} \\%
U \arrow[swap]{u}{\pi}  \arrow{r}{\tau_\beta} & U 
\end{tikzcd} 
\]
where 
\begin{align}
& \tau_\beta(l,k,j) = (l',k',j') 
\label{taubeta} \\
& k' = k-j-l \nonumber \\
& j' = k \nonumber \\
& l' = \begin{cases}
0,\ \text{ if } j'=k' = 0, \\
1+\lfloor k'\alpha -j'\frac{\alpha}{\beta} \rfloor  , \text{ otherwise.}
\end{cases} \nonumber 
\end{align}
Note that $f(l',k',j')\in [0,1)$, and is zero only for $(k',j')=(0,0)$. 
It is clear that $\tau_\beta(l,k,j)= G(l,k,j) - \lfloor f\big(G(l,k,j)\big) \rfloor (1,0,0)$, but the expression in \eqref{taubeta} will benefit future discussions. 
Finally, note that if $(k',j')=0$, i.e., $k=0$ and $-j-l=0$, then
\[ f(l,0,j)= l +j\frac{\alpha}{\beta} \in [0,1) \]
implies that if $j\neq 0$, then $j=-1$ and $l=1$. This means that the final symbol of the $\beta$-expansion, for any $\omega_n(l,k,j)\in \finb$, is precisely $n=\lfloor f\big(G(1,0,-1)\big)\rfloor$. 

On the known studies of shift radix systems \cite{ak_all,taka1} the use of maps in $\Z^{d-1}$, where $d$ is the degree of minimal polynomial of $\beta$, is analogous to the map $\tau_\beta$. The matrix $G$ is similar to the companion matrix of the minimal polynomial, and was obtained from the parametrization of the $\beta$-fractional part of {\em large} natural numbers in Corollary \ref{omega}. 

\begin{proposition}
If $\omega_n(l,k,j)\in \finb$, then its $\beta$-expansion has the suffix $10n$.
\label{suffix1}
\end{proposition}

\begin{proof}
We just argued that the only pre-image of $(0,0,0)$ in $U^*$ is $(1,0,-1)$.
As we will go backwards, let us write $v_3=(1,0,-1)$. Write 
$$H=\begin{pmatrix}
1 & 0  & 1 \\ 0 & 0 & n \\ -1 & -n & n-1 \end{pmatrix}$$ 
for the adjoint of $G$. 

In order to get a pre-image of $v_i=(l_i,k_i,j_i)\in U^*$, let $0\leq r_{i-1}<n$ 
be the residue of the division $l_i+j_i$ by $n$. If $f(v_i)\geq \alpha$ and 
$r_{i-1}=0$, set $s_{i-1}=0$. Otherwise set $s_{i-1}=n-r_{i-1}$. 

Define 
\[ v_{i-1}= \frac{H}n \begin{pmatrix} l_i+s_{i-1} \\ k_i \\ j_i \end{pmatrix}= 
\begin{pmatrix} \frac{l_i+s_{i-1}+j_i}{n} \\
j_i \\ 
j_i-k_i - \frac{l_i+s_{i-1}+j_i}{n} \end{pmatrix} \ ,\]
then $f(v_{i-1})= \frac 1{\beta} (f(v_i)+s_{i-1})$, since $fH= \frac n{\beta} f$. 

Note that $s_{i-1}$ is the integer part of $\beta f(v_{i-1})$, which makes the 
choice of the pre-image unique in $U$. For $v_3=(1,0,-1)$, we have $r_2=0$ and since $f(v_3)=1-\frac{\alpha}{\beta}>\alpha$, we get $s_2=0$.
Moreover $v_2 =\frac 1n H\begin{pmatrix} 1 \\ 0 \\ -1 \end{pmatrix} = \begin{pmatrix} 0 \\ -1 \\ -1 \end{pmatrix}$. We have $r_1=n-1$ and $s_1=1$. 
\end{proof} 

The argument in the proof of Proposition \ref{suffix1} depends on the choice of inverse images for points in $U$. Remarkably, the map $\tau_\beta$ is at most 2-to-1 in $U$, whereas $T_\beta : [0,1)\to [0,1)$ is at least $(n-1)$-to-1. 

\begin{proposition}
Let $(l_0,k_0,j_0),\ (l_1,k_1,j_1)\in U$ be such that $\tau_\beta(l_0,k_0,j_0)=(l_1,k_1,j_1)$. If $l_1+j_1\not\equiv 0\pmod n$, then $\tau_\beta$ is 1-to-1. If $l_1+j_1\equiv 0 \mod n$ and $f(l_0,k_0,j_0)\in[\alpha/\beta,1-\alpha/\beta)$, then $\tau_\beta$ is 1-to-1. If $l_1+j_1\equiv 0 \mod n$ and $f(l_0,k_0,j_0)<\alpha/\beta$, then $\tau_\beta$ is 2-to-1, $(l_0+1,k_0,j_0-1)\in U$ and $\tau_\beta(l_0+1,k_0,j_0-1) = (l_1,k_1,j_1)$. 
If $l_1+j_1\equiv 0 \mod n$ and $f(l_0,k_0,j_0)> 1-\alpha/\beta$, then $\tau_\beta$ is 2-to-1, $(l_0-1,k_0,j_0+1)\in U$ and $\tau_\beta(l_0-1,k_0,j_0+1) = (l_1,k_1,j_1)$. 
\end{proposition}

The proof is clear and we omit it. The {\em exceptional} non-injectiveness of $\tau_\beta$ in the domain $U$ is explained by the conditions 
\begin{align*}
f(l_0,k_0,j_0) &  <\alpha/\beta \ ,\\
f(1,0,-1) & = 1-\alpha/\beta 
\end{align*}
together with the fact that $(1,0,-1)$ is the only pre-image of the origin in $U^*$. 

We can point out another important characteristic of the points of non-injectiveness of $\tau_\beta$. Since
\[ \lfloor f\big( G(l_0+1,k_0,j_0-1)\big)\rfloor= \lfloor f(nl_0+n-k_0, k_0-j_0-l_0,k_0) \rfloor = n + \lfloor f(l_1,k_1,j_1) \rfloor \ , \]
the letter for the point parametrized by $(l_0+1,k_0,j_0-1)$ is $n$, while the letter for the point parametrized by $(l_0,k_0,j_0)$ is $0$. Consequently

\begin{proposition}
Let $(l,k,j)\in U^*$ and let its orbit $\big(\tau_\beta^p(l,k,j)\big)_{p\geq 0}$ be coded by a sequence $(s_p)$, finite or infinite, such that for each $p$ such that $s_p=0$, the replacement by $s_p=n$ generates a forbidden word. Then $\tau_\beta$ is invertible along the forward orbit of $(l,k,j)$.
\end{proposition}

We turn now to the study of the orbits $U^*$ for $n\geq 4$. 

\begin{lemma}
Let $(l,k,j)\in U^*$, $n\geq 4$, $m=\max\{|k|,|j|\}$ and $\mu=\min\{|k|,|j|\}$. Then 
\[ -\frac{\mu}{n-2}-\frac{m-\mu}{n-1} <l< 1+ \frac{\mu}{n-2}+ \frac{m-\mu}{n-1} \ .\] 
\label{lsmall}
\end{lemma}

\begin{proof}
By its formula, $l$ is largest when $k>0>j$, and $m=k$, in which case
\begin{align*}
l & = 1+\lfloor k\alpha-j\alpha/\beta\rfloor \leq 1 + \lfloor \mu \alpha(1+1/\beta)+(m-\mu)\alpha \rfloor \\ 
& < 1+ \frac{\mu}{n-2}+\frac{m-\mu}{n-1} \ ,
\end{align*}
where we have used that $\alpha$ is monotonically decreasing as a function of $n$, 
\[ \alpha(1+1/\beta)< \frac 1{n-2},\ \text{ and } (n-1)\alpha <1 \ . \]
On the other hand, $l$ is smallest when $k<0<j$, $|k|=m$ and $j=\mu$:
\begin{align*}
l &\geq 1+\lfloor -\mu(\alpha+\alpha/\beta)-(m-\mu) \alpha \rfloor \\
& > -\mu\alpha(1+1/\beta)-(m-\mu)\alpha \\ 
& > -\frac{\mu}{n-2}- \frac{m-\mu}{n-1}
\end{align*}
\end{proof}

\begin{remark}
Let us collect here the study of $l$ for small $m$.
If $kj>0$ and $m\leq n$, then for $k>0$
\begin{align*}
l= 1+\lfloor k\alpha - j\alpha/\beta \rfloor \leq 1+\lfloor n(\alpha-\alpha/\beta)\rfloor = 1 \ , 
\end{align*}
and similarly $l=0$ for $k<0$. 

If $m\leq n-2$ and $kj<0$, then $l=1$ for $k>0$ and $l=0$ for $k<0$. 
If $k=0$ and $0<j\leq n^2-2$, then $l=0$. If $k=0$ and $2-n^2\leq j<0$, then $l=1$. 
If $j=0$ and $1-n\leq k<0$, then $l=0$. If $j=0$ and $0<k\leq n-1$, then $l=1$.

In particular, if $m=1$, then $0\leq l\leq 1$. 
\label{obs_studyl}
\end{remark}

\begin{corollary}
Let $(l,k,j)\in U^*$, $n\geq 4$ and $m=\max\{|k|,|j|\}$. If $m\geq 2$, then $m>|l|$. 
\label{cor1}
\end{corollary}

\begin{proof}
Suppose $m\leq |l|$, then 
\[ m < 1+ \frac{\mu}{n-2}+\frac{m-\mu}{n-1} \leq 1+\frac{m}{n-2} \ . \]
We conclude that $(m-1)(n-3)<1$, which is false.
\end{proof}

Remark \ref{obs_studyl} and Corollary \ref{cor1} show that the norm in $U$: $\max\{|l|,|k|,|j|\}=\max\{|k|,|j|\}=m$.

\begin{lemma}
\label{first_steps}

\begin{enumerate}
\item[a)] If $k>0$ and $j<0$, then 
\[ k-j-l > m\frac{n-2}{n-1} + \mu \frac{n^2-3n+1}{n^2-3n+2}-1 >0\ .\]
\item[b)] If $k<0$ and $j>0$, then 
\[ k-j-l< -m \frac{n-2}{n-1} - \mu \frac{n^2-3n+1}{n^2-3n+2}  \ .\]
\end{enumerate}
\end{lemma}

\begin{proof}
To see (a), write $\mu=\min\{k,-j\}$, and note that $l\geq 1$. Then
\begin{align*}
k-j-l & = |k|+|j|-|l|> m+\mu-1-\frac{\mu}{n-2}-\frac{m-\mu}{n-1} \\
k-j-l & > m\frac{n-2}{n-1} + \mu \frac{n^2-3n+1}{n^2-3n+2}-1>0 \ . 
\end{align*}
To see (b), note that $l\leq 0$, and
\begin{align*}
k-j-l & = -|k|-|j|+|l| < -(m+\mu)+\frac{\mu}{n-2}+\frac{(m-\mu)}{n-1}  \\
k-j-l& < -m \frac{n-2}{n-1} - \mu \frac{n^2-3n+1}{n^2-3n+2} \ .
\end{align*}
\end{proof}

\begin{proposition}
\label{fivesteps}
Let $(l_i,k_i,j_i)$, $0\leq i \leq 4$, $(l_0,k_0,j_0)\in U^*$ and 
$\tau_\beta(l_i,k_i,j_i)= (l_{i+1},k_{i+1},j_{i+1}) \in U^*$ be five points in the orbit of $(l_0,k_0,j_0)$. If $n\geq 4$, $k_0>0$, $j_0>0$ and $k_1<0$, then $k_2<0$, $k_4>0$ and $k_5>0$. 
\end{proposition}

\begin{proof}
Note that $j_{i+1}=k_i$, which implies that $j_1>0$, $j_2<0$, $j_3<0$, $j_5>0$ and $j_6>0$. We keep the notation $\mu_i = \min\{|k_i|,|j_i|\}$ and 
$m_i =\max\{|k_i|,|j_i|\}$. 

Since $k_1<0$ and $j_1>0$, then $l_1\leq 0$ and Lemma \ref{first_steps}(b) can be applied 
\[ k_2=k_1-j_1-l_1 < -m_1\frac{n-2}{n-1} - \mu_1 \frac{n^2-3n+1}{n^2-3n+2}\ \]
proving the first statement. Note that $j_2=k_1<0$, and 
\[ \frac 1{\beta} < \frac{n-2}{n-1} < \frac{n^2-3n+1}{n^2-3n+2} \ . \]
Since either $k_1=-m_1$ or $k_1=-\mu_1$, we have $k_2-j_2/\beta<0$, implying that $l_2\leq 0$. 

From this, one gets
\[ k_2-j_2 \leq  k_2-j_2-l_2 = k_3 \ . \]

We now split the proof in three cases. 

\noindent {\it Case i. $k_3<0$}

The inequality $k_2-j_2\leq k_3$ implies $|k_3|\leq |k_2-j_2|$. On the other hand, $k_2-j_2=-j_1-l_1=-k_0-l_1$. Thus $k_2-j_2\leq k_3<0$ implies $l_1>-k_0$, and since we know that $l_1\leq 0$, we get $-k_0<l_1\leq 0$. Finally $k_2-j_2<0$ also implies
\[ |k_3|\leq |k_2-j_2| < |k_2|=|j_3| \ , \]
meaning that $m_3=|j_3$ and $\mu_3=-k_3$. 

Consider $m_2=\max\{|j_2|,|k_2|\}$. If $m_2\leq n$, then $m_3=|j_3|=|k_2|\leq n$. 
Then Remark \ref{obs_studyl} yields 
$l_3= 0$ and finally $k_4=k_3-j_3>0$. 

Suppose $m_2>n$. If $l_3\leq 0$, $k_4=k_3-j_3+(-l_3)>0$. 

If $l_3>0$, then 
$l_3<1+ \frac{\mu_3}{n-2}+ \frac{m_3-\mu_3}{n-1}$ from Lemma \ref{lsmall}, so that 
\begin{align*}
k_4=k_3-j_3-l_3 &> -\mu_3+m_3-1-\frac{m_3}{n-1}-\frac{\mu_3}{n^2-3n+2}
\\
k_4 &> m_3 \frac{n-2}{n-1} -\mu_3 \frac{n^2-3n+3}{n^2-3n+2}-1 \ . 
\end{align*}
On the other hand, 
\[ l_3 = 1+\lfloor k_3\alpha - j_3\alpha/\beta\rfloor >0 \; \Rightarrow \; k_3>\frac{j_3}{\beta} \ ,\]
which yields $m_3=|j_3|>\beta |k_3|=\beta\mu_3>n\mu_3$. Then
\[ k_4> \mu_3 \frac{n(n-2)^2-n^2+3n-3}{n^2-3n+2}-1=\mu_3\frac{n^2-4n+3}{n-2}-1>0 \ , \]
as $n\geq 4$. 

Applying Lemma \ref{first_steps}(a), we get $k_5>0$ and $j_5=k_4>0$. 

\noindent {\it Case ii. $k_3=0$} 

Then $j_3=k_2<0$ and $l_3\geq 1$, but
\[ l_3<1-k_2 \frac{\alpha}{\beta} < 1-\frac{k_2}{n^2-n} \ . \]
Then 
\[ k_4=-j_3-l_3> -k_2 \frac{n^2-n+1}{n^2-n}-1 > 0 \ ,\]
for $k_2\leq -2$. If $k_2=-1$, then $(l_3,k_3,j_3)=(1,0,-1)$ would imply $(l_4,k_4,j_4)=(0,0,0)\not\in U^*$, contrary to the hypothesis.
From the point $(l_4,k_4,0)$, it is also clear that $k_5=k_4-l_4>0$ and $j_5=k_4>0$.

\noindent{\it Case iii. $k_3>0$.} 

Since $j_3<0$, then $k_4>0$ follows from Lemma \ref{first_steps}(a). 
We also have $|j_3|=m_3$ and $\mu_3 = k_3$. Lemma \ref{first_steps}(a) explicitly gives
\[ k_3-j_3-l_3 > m_3\frac{n-2}{n-1}+\mu_3 \frac{n^2-3n+1}{n^2-3n+2} -1 \ , \]
which implies, when $|j_3|>k_3$
\[ m_3-l_3> m_3\frac{n-2}{n-1}-\mu_3 \frac{ 1}{n^2-3n+2}-1 \ . \]

Suppose that $k_5< 0$:
\begin{align*}
& k_4-j_4-l_4<0 \\
& k_4-j_4<l_4=1+\lfloor k_4 \alpha- j_4\alpha/\beta\rfloor \\
& k_4-j_4 -\lfloor k_4 \alpha- j_4\alpha/\beta\rfloor <1 \\
& k_4-j_4 -\lfloor k_4 \alpha- j_4\alpha/\beta\rfloor \leq 0 \\
& k_4-j_4 \leq \lfloor k_4 \alpha- j_4\alpha/\beta\rfloor <  k_4 \alpha- j_4\alpha/\beta \\
& k_4(1-\alpha)< j_4 (1-\alpha/\beta) \\
& (m_3+\mu_3-l_3)(1-\alpha) < \mu_3(1-\alpha/\beta) \\
& (m_3-l_3)(1-\alpha)< \mu_3\alpha(1-1/\beta) \\
& 1-\alpha > \frac{n-2}{n-1}\; \ , \alpha(1-1/\beta) < \frac{n}{n^2-1} \\
& (m_3-l_3)(n-2)< \mu_3 \frac{n}{n+1} < \mu_3 
\end{align*}
We eventually get 
\[ m_3 < \mu_3 \frac{n}{(n-2)^2} +\frac{n-1}{n-2} \ , \]
which is clearly false.

If $k_5=0$,  a similar calculation entails the estimate
\[ m_3 < \frac{2-\alpha}{1-\alpha-\alpha/\beta} + \mu_3 \frac{2\alpha-\alpha/\beta}{1-\alpha-\alpha/\beta} \ ,\]
leaving a handful of possibilities for the pair $(k_3,j_3)$ each of which is checked contradictory.
\end{proof}

The situation described in Proposition \ref{fivesteps} entails that $\tau_\beta^5$ or $\tau_\beta^6$ is an endomorphism of the positive cone in the plane $\pi(l,k,j)=(k,j)$ into itself. 

A pre-periodic orbit must have a point where $\tau_\beta$ is not one-to-one, which forces it to visit at least one of the intervals $(0,\alpha/\beta)$ or $(1-\alpha/\beta,1)$. If $\|\cdot \|$ denotes the distance to the closest integer, we may say that a pre-periodic orbit has a point $(l,k,j)\in U^*$ such that $0<\|f(l,k,j)\|<\alpha/\beta$. 

Recall that $1,\ \alpha$ and $\alpha/\beta$ are rationally independent. We now use W. Schmidt's study on simultaneous approximations of two rationally independent algebraic number by rationals.
From Corollary 1 of \cite{schimdt67}, there is a finite number of triples $(l,k,j)\in U^*$ such that 
\begin{equation}
\label{cond1}
 0<\|f(l,k,j)\|<\frac{\alpha}{\beta}< \frac 1{n^{2+\epsilon}} \ , 
\end{equation}
with $\max\{|k|,|j|\}\leq n$. In fact, a direct solution of the inequalities $0<\|f(0,k,j)\|<\frac{\alpha}{\beta}$ gives two triples
$v_0=(0,-1,-n)$ and $v_1=(0,1,n+1)$. These will generate two solutions for $l=1$, 
\begin{align*}
& u_0=(1,0,0)-(0,-1,-n) = (1,1,n) \ , \\
& u_1=(1,0,0)-(0,1,n+1)=(1,-1,-n-1) 
\end{align*}
These four solutions in the ball contained in $U^*$ satisfying $\max\{|k|,|j|\}\leq n+1$ and \eqref{cond1} are on the pre-orbit, under $\tau_\beta$ of $(1,n-1,-1)$, which projects to $\alpha^2$, and indeed $f(u_1)$ and $f(v_0)$ are pre-images of $\alpha^2$. 

We are now close to finish the argument proving a contraction condition on the map $\tau_\beta$. Before that, we single out some points in $U^*$ which belong to periodic orbits. They have the form $(1,k,k+j-1)$, with $k,j>0$ odd, $k<n-1$ and $k+j-1<n$. In the spirit of section \ref{cubic}, we can show that the $\beta$-expansion of $f(1,k,k+j-1)$ will be periodic. For some word $w_n(k,j)$  
\[ [f(1,k,k+j-1)]_\beta = w_n(k,j) [f(1,k,k+j-1)]_\beta \ . \]
Indeed, letting $u_n(k,j)= (n-k)j(k+j)(k-1)(n-j-2)(n-k-j)$, then
\begin{align*}
& w_n(k,j) = u_n(k,j) u_n(k-2,j+2)\ldots u_n(3,j+k-3) (n-1)(j+k-1)\\ 
& \;\; (j+k)0(n-j-k) u_n(j+k-1,1)u_n(j+k-3,3)\ldots u_n(k-2,j+2) \ . 
\end{align*}
Note that the length of $w_n(k,j)$ is congruent to $5 \mod 6$.

The following observation is useful.

\begin{lemma}
$\beta^{-1}\not\in\Z[\beta]$. 
\label{beta_1}
\end{lemma}

\begin{proof}
$\beta$ is not a unity, and
\[ \beta^{-1} = \frac{\beta^2}{n} - \frac{n+1}{n} \beta+1 \ . \]
If $\beta^{-1}\in \Z[\beta]$, $\beta$ would be root of a polynomial of degree 2 with coefficients in $\Z$, which is an absurd.
\end{proof}

\begin{proposition}
The numbers whose $\beta$-fractional parts are of the form $f(1,k,k+j-1)$, for $0<k\leq n$ and $0<k+j-1\leq n$, are not positive integers.
\label{periodic_orb}
\end{proposition}

\begin{proof}
Rewriting eq.\eqref{lema1}:
\[ kn = (k-1)\beta+(n-1)+f(1,k-1,0) \ . \]
Therefore
\[ kn-\alpha+(k+j-1)\frac{\alpha}{\beta} \]
has $\beta$-expansion given by $(k-1)(n-1).w$, where $w$ is the $\beta$-expansion of $f(1,k,k+j-1)$. 

If there existed $N\in \N$ with the same $\beta$-fractional part as $kn-\alpha+(k+j-1) \frac{\alpha}{\beta}$, then
\[ N-kn+\alpha-(k+j-1)\frac{\alpha}{\beta} \in \Z[\beta] \ . \]
Since $\alpha=\beta-n$, this contradicts Lemma \ref{beta_1}. 
\end{proof}

On the other hand, from a point of the form $(1,k,k+\delta)$, with $\delta\geq 0$ odd or $k$ even, similar unfolding calculations for 
$[f(1,k,k+\delta)]_\beta$ inform that these points have finite orbits. From Proposition \ref{fivesteps}, given any point 
$(l,k,j)$ in the ball $\max\{|k|,|j|\}\leq n$, its orbit will visit the first quadrant 
$0\leq k,j\leq n$. Finally, from Remark \ref{obs_studyl}, if $(l,k,j)$ is such that 
$0< j<k\leq n$, then $l=1$ and $j_1=k\geq k_1=k-j+1$.

\begin{proposition}
There is a ball $B_p=\max\{|k|,|j|\}\leq p$ which is invariant under $\tau_\beta$. Every point outside $B_p$ is eventually mapped into $B_p$.  
\label{contrai}
\end{proposition}

\begin{proof} 
We write $\tau_\beta$ in terms of the following basis in $\R^3$:
\[ \gamma=\big\{(1,-\alpha,\alpha{\beta}^{-1}), (1,\alpha^{-1},0), (1,0,-\beta{\alpha}^{-1})\big\} \ . \]
If $(l,k,j)= \lambda_1 (1,-\alpha,\alpha{\beta}^{-1})+\lambda_2(1,\alpha^{-1},0) + \lambda_3(1,0,-\alpha^{-1}\beta)$, then 
we see that $\lambda_1= \frac{\beta}{\beta+\alpha+\alpha^2} f(l,k,j)\in (0,1)$ for any $(l,k,j)\in U^*$ and
\[ P= \begin{pmatrix} 
1 & 1 & 1 \\
-\alpha & \alpha^{-1} & 0 \\
\alpha\beta^{-1} & 0 & -\beta\alpha^{-1} \end{pmatrix} \ \]
changes basis so that
\[ G_\gamma = \begin{pmatrix} \beta & 0 & 0 \\ -2\alpha^2\beta^{-1} & 1-\alpha & n \\
2\alpha^2\beta^{-1} & - \beta^{-1} & 0 \end{pmatrix} \ . \]
Let us introduce the abbreviation $\xi=(\beta+\alpha+\alpha^2)$, then $\tau'_\beta \stackrel{\rm def}{=} P^{-1} \tau_\beta P$ is given by
\[ \tau'_\beta(\lambda_1,\lambda_2,\lambda_3) =
\begin{pmatrix}
\beta\lambda_1 - \frac{\beta \lfloor \lambda_1 \xi \rfloor}{\xi} \\
-2\frac{\alpha^2}{\beta}\lambda_1 -\frac{\alpha^2\beta \lfloor \lambda_1 \xi \rfloor }{\xi} +(1-\alpha)\lambda_2 +n\lambda_3 \\
2\frac{\alpha^2}{\beta}\lambda_1-\frac{\alpha^2 \lfloor \lambda_1 \xi \rfloor }{\beta\xi}  -\frac{\lambda_2}{\beta} 
\end{pmatrix} \ . \]

We notice that the subspace $S$ spanned by $\{(1,\alpha^{-1},0),(1,0,-\beta\alpha^{-1})\}$ is quasi-invariant under $\tau_\beta'$. Restricted to it, $\tau_\beta'$ is a perturbation of a linear map, which is a contraction
represented by the $2\times 2$ matrix: $A=\begin{pmatrix} (1-\alpha) & n \\
-\beta^{-1} & 0 \end{pmatrix}$. This implies that there is an open set containing the origin which is invariant under $\tau_\beta'$. Furthermore, on the complement of this set, $\tau_\beta'$ is a contraction on the subspace $S$. This observation proves the second statement.

A simple perturbative calculation  yields that
\[ B_\gamma =\{(\lambda_1,\lambda_2, \lambda_3)\ :\ 0<\lambda_1<1, |\lambda_2|\leq n, |\lambda_3|\leq 1 \} \cap (P^{-1} U) \ \]
is invariant. Indeed, if $\Delta= \begin{pmatrix} -2\frac{\alpha^2}{\beta}\lambda_1 -\frac{\alpha^2\beta \lfloor \lambda_1\xi\rfloor}{\xi} \\
2\frac{\alpha^2}{\beta}\lambda_1 - \frac{\alpha^2\lfloor \lambda_1\xi\rfloor}{\beta\xi} 
\end{pmatrix}$, the solution to 
\[ \Delta+A\begin{pmatrix} \lambda_2 \\ \lambda_3 \end{pmatrix} = \begin{pmatrix} n \\ 1 \end{pmatrix} \]
satisfies $\lambda_2<-\beta$ and $\lambda_3>1$. 

Finally $B_n\subset P B_\gamma$, because $|\lambda_3|\approx \alpha\beta^{-1} |j|$ and $|\lambda_2|\approx \alpha |k|$. $|k|\leq n$ is sufficient for $|\lambda_2|\leq n$ and $|j|\leq n$ implies $|\lambda_3|\leq 1$. 
\end{proof}

\begin{proof}[Proof of Theorem \ref{um}]
From Proposition \ref{fivesteps}, a point in a periodic orbit will eventually visit the first quadrant. In the basis $\gamma$, the coordinates satisfy
\[ \lambda_1+\lambda_2+\lambda_3 = l \ , \; 0<\lambda_1<1  \ . \]
Since eventually the orbit visits the strip $l=1$, we obtain 
\[ 0<\lambda_2<1 \ , \; -1<\lambda_3<0 \ ,\]
since $\lambda_2\approx \alpha k$ and $\lambda_3\approx - \alpha j/\beta$. Hence every periodic orbit visits $B_n$. 

If $(l,k,j)$ is a point of a periodic orbit such that $0\leq j<k\leq n$ then $l=1$ and $j_1=k>k_1=k-j+l$ and $l_1=1$, from Remark \ref{obs_studyl}. Every periodic orbit having a point of the form $(1,k,k+\delta)$, with $\delta\geq 0$ does not correspond to the $\beta$-fractional of a natural number, from Proposition \ref{periodic_orb}.
The possible periodic orbit remaining is that corresponding to $\alpha^2$, whose coordinates are $(1,n-1,-1)$. We have shown, however, that $\alpha^2\in \finb$ for every $n\equiv 1 \mod 3$. 
\end{proof} 

\section{Final remarks}
\label{ultima}

In this section, we try to gain some perspective for future work.

The map $f: \Z^3\to \R$, which helped us to obtain a conjugate map to the $\beta$-transformation was unveiled through detailed calculations of the $\beta$-fractional parts of the sucessor of a natural number for the cubic family 
\eqref{pol1}. One question is whether this approach generalizes to other cubic Pisot polynomials or even of Pisot polynomials of higher degree. Proposition \ref{prep_induction} plays an analogous rôle as Corollary 2.14 from \cite{taka1}. 

Given $\beta>1$ a cubic Pisot number, it is rather natural to use the generalization of $T_\beta: \pi(U)\to \pi(U)$ where $U$ is given in \eqref{pontosdez3} and the functional $f$ has coefficients corresponding to the left eigenvector of the companion matrix of the minimal polynomial $X^3-a_2X^2-a_1X-a_0$
for $\beta$:
\[ \begin{bmatrix} 0 & 1 & 0 \\ 0 & 0 & 1 \\ a_0 & a_1 & a_2 \end{bmatrix} \ , \]
that is, $f=[r_1,r_2,1]$, where $\beta r_1=a_0$, $\beta r_2=r_1+a_1$ and $a_2+r_2=\beta$. 

If $1$ has a finite $\beta$-expansion of length $k$, write $n=\lfloor \beta \rfloor$ and $\alpha=\beta-n$. The $\beta$-expansion of $\alpha$ is the proper suffix of $[1]_{\beta e}$ of length $k-1$. 
A direct calculation shows that
\begin{align*}
& \alpha= r_2+(a_2-n) \text{ and } \\
& \alpha^2= r_1+r_2(a_2-2n)+(a_1+(a_2-n)^2) \ ,
\end{align*}
which means that $\alpha^2 \in \Z[\beta]$, and that $\{\alpha,\alpha^2\}$ is a linear independent set in $\Z[\beta]$ as a $\Z$-module. Therefore $\{1,\alpha,\alpha^2\}$ is a basis for $\Z[\beta]$. 
Note that $\alpha=f(0,1,a_2-n)$.

If the union of pre-images of the origin 
\[ F_\beta =\left\{ (l,k,j)\in \Z^3\ |\ \exists k>0,\ T_\beta(l,k,j)=0 \right\} \]
contains the integer linear combinations of $T_\beta^{n-1}(0,1,a_2-n)$, for any $n\in \N$ (again cf. Proposition 3.2 from \cite{taka1}), then $(F_1)$ holds. 

However, the minimal polynomial $x^3-x^2-3x-2$, which was proved to be a Pisot polynomial by Bassino \cite{bassino}, provides the following finite $\beta$-expansions:
\begin{align*}
[1]_{\beta e} & = 21012 \\
[\alpha]_{\beta e} & = 1012 \\
[\alpha^2]_{\beta e} & = 0 11111011001012
\end{align*}
and yet $7\notin \finb$. We observe that the $\beta$-fractional part in the $\beta$-expansion of 7 is
$2-2\alpha-2\tau_\beta(\alpha)$.

In terms of the set of periodic orbits in $U^*\subset \Z^d$, say $P_\beta$, lack of injectivity of $T_\beta$ is necessary when 
$P_\beta$ is not invariant under $T_\beta$, that is condition (i) in Theorema 3.3 does not hold. Injectivity, and lack thereof, was mentioned in Proposition 3.2 and 3.3 of this paper. 

We have to leave open the question of a general alternative formulation of sufficient conditions for $(F_1)$.









\end{document}